\documentclass[11pt,leqno]{article}
\usepackage{a4}
\usepackage[T1]{fontenc}
\usepackage[english]{babel}
\usepackage{latexsym}
\usepackage{amsmath}
\usepackage{amssymb}
\usepackage{mathrsfs}
\usepackage{fig4tex}
\usepackage{color}
\usepackage{cite}
\usepackage{titling}
\usepackage{bbold}
\makeatletter
\@addtoreset{equation}{section} 
\makeatother
\definecolor{move}{rgb}{.3,.1,.8}
\newcommand{\ud}{\mathrm{d}}
\newcommand{\vfi}{\varphi}

\newcommand{\re}{\mathrm{Re}}
\newcommand{\im}{\mathrm{Im}}
\newcommand{\e}{\mathrm{e}}
\newcommand{\id}{\mathrm{Id}}
\newcommand{\Rn}{\mathbb{R}^{n}}

\newcommand{\p}{\partial}
\newcommand{\pn}{\partial_{\nu}}
\newcommand{\Ci}{\mathscr{C}^{\infty}}
\newcommand{\Cc}{\mathscr{C}_{c}^{\infty}}

\newenvironment{pr}{\vspace{5pt}\textbf{{\small Proof :}}\\}{\hspace{\stretch{1}}\rule{1ex}{1ex}\vspace{5pt}}
\newtheorem{thm}{Theorem}[section]
\newtheorem{pro}{Proposition}[section]

\newtheorem{lem}{Lemma}[section]
\newtheorem{rem}{Remarks}[section]

\usepackage{fancyhdr}
\fancypagestyle{plain}{\fancyhead{}}
\title{Weak stabilization of a transmission Euler-Bernoulli plate equation with force and moment feedback}
\author{{FATHI HASSINE}\\ \textit{ UR Analysis and Control of PDE (05/UR/15-01)}\\ \textit{D\'epartement de Math\'ematiques, Facult\'e des sciences de Monastir,}\\ \textit{5019 Monastir, Tunisie}\\ \textit{email:} \texttt{fathi.hassine@fsm.rnu.tn}}
\date{}
\begin{document}
\maketitle
\begin{center}
\abstract{In this paper we will study the asymptotic behaviour of the energy decay of a transmission plate equation with force and moment feedback. Precisly, we shall prove that the energy decay at least logarithmically over the time. The method consist to use the classical second order Carleman estimate to estabish a resolvent estimate which provide by the famous Burq's result~\cite{Bur} the kind of decay above mentionned.}
\end{center}
\textbf{Key words and phrases: }Transmission problem, boundary stabilization, Euler-Bernoulli plate equation, energy decay, Carleman estimates.
\\
\textbf{Mathematics Subject Classification:} \textit{35A01, 35A02, 35M33, 93D20}.
\tableofcontents
\newpage
\section{Introduction and statement of results}
\begin{figure}[htbp]
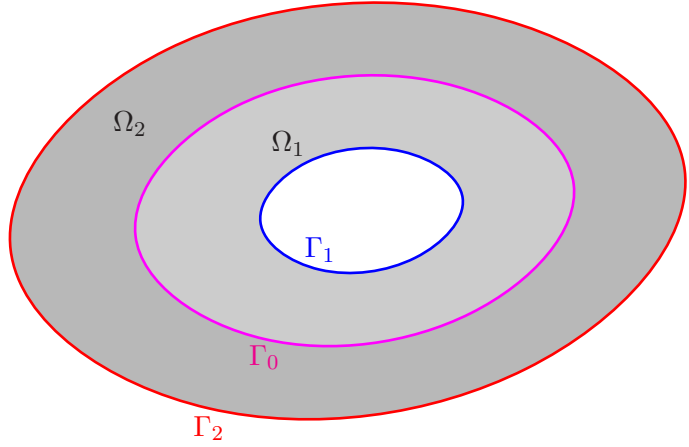

\figinit{1.5pt}
\figpt 0:(20,0)
\figpt 1:(-70,-10)
\figpt 2:(-20,-50)
\figpt 3:(60,-40)
\figpt 4:(100,10)
\figptsym 5:=3/1,4/
\figptsym 6:=2/1,4/
\figpthom 7:=0/1,0.35/
\figpthom 8:=0/2,0.35/
\figpthom 9:=0/3,0.35/
\figpthom 10:=0/4,0.35/
\figpthom 11:=0/5,0.35/
\figpthom 12:=0/6,0.35/
\figpthom 14:=0/1,0.7/
\figpthom 15:=0/2,0.7/
\figpthom 16:=0/3,0.7/
\figpthom 17:=0/4,0.7/
\figpthom 18:=0/5,0.7/
\figpthom 19:=0/6,0.7/
\figpt 13:(-40,25)
\figpt 20:(0,20)
\psbeginfig{}
\psset(width=1)
\psset(fillmode=yes,color=0.72)
\pscurve  [1,2,3,4,5,6,1,2,3]
\psset(fillmode=yes,color=0.8)
\pscurve  [7,8,9,10,11,12,7,8,9]
\psset(fillmode=yes,color=1)
\pscurve  [14,15,16,17,18,19,14,15,16]
\psset(fillmode=no)
\psset(color=\Magentargb)
\pscurve  [7,8,9,10,11,12,7,8,9]
\psset(color=\Redrgb)
\pscurve[1,2,3,4,5,6,1,2,3]
\psset(color=\Bluergb)
\pscurve  [14,15,16,17,18,19,14,15,16]
\psendfig
\figvisu{\figBoxA}{}{
\figwrites 20:$\Omega_{1}$(0.15)
\figwrites 13:$\Omega_{2}$(0.15)
\figwrites 8:\textcolor{magenta}{$\Gamma_{0}$}(1.5)
\figwriten 15:\textcolor{blue}{$\Gamma_{1}$}(2)
\figwrites 2:\textcolor{red}{$\Gamma_{2}$}(2)
}
\centerline{\box\figBoxA}
\centerline{\box\figBoxA}
\caption{Geometrical situation of the transmission problem.}
\label{fig1}
\end{figure}

Let $\Omega\subset\Rn$ be an open, bounded connected domain with smooth boundary $\Gamma=\Gamma_{1}\cup\Gamma_{2}$ where $\Gamma_{1}$ and $\Gamma_{2}$ are two non empty component of $\Gamma$ such that $\Gamma_{1}\cap\Gamma_{2}=\emptyset$.
\\
Let $\Omega_{1}\subset\Omega$ be an open domain with smooth boundary $\p\Omega_{1}=\Gamma_{1}\cup\Gamma_{0}$ where $\Gamma_{1}\cap\Gamma_{0}=\Gamma_{2}\cap\Gamma_{0}=\emptyset$. Then $\Omega_{2}=\Omega\backslash\overline{\Omega}_{1}$ is an open connected domain with boundary $\p\Omega_{2}=\Gamma_{2}\cup\Gamma_{0}$ (See Figure~\ref{fig1}).
\\
We are going to study the following mixed boundary value problem
\begin{equation}\label{1}
\left\{
\begin{array}{lll}
\p_{t}^{2}u_{1}+c_{1}^{2}\Delta^{2}u_{1}=0&\textrm{in}&\Omega_{1}\times]0,+\infty[,
\\
\p_{t}^{2}u_{2}+c_{2}^{2}\Delta^{2}u_{2}=0&\textrm{in}&\Omega_{2}\times]0,+\infty[,
\\
\textcolor{magenta}{u_{1}=u_{2}}&\textcolor{magenta}{\textrm{on}}&\textcolor{magenta}{\Gamma_{0}\times]0,+\infty[},
\\
\textcolor{magenta}{\pn u_{1}=\pn u_{2}}&\textcolor{magenta}{\textrm{on}}&\textcolor{magenta}{\Gamma_{0}\times]0,+\infty[},
\\
\textcolor{magenta}{c_{1}\Delta u_{1}=c_{2}\Delta u_{2}}&\textcolor{magenta}{\textrm{on}}&\textcolor{magenta}{\Gamma_{0}\times]0,+\infty[},
\\
\textcolor{magenta}{c_{1}\pn \Delta u_{1}=c_{2}\pn \Delta u_{2}}&\textcolor{magenta}{\textrm{on}}&\textcolor{magenta}{\Gamma_{0}\times]0,+\infty[},
\\
\textcolor{blue}{\Delta u_{1}=0}&\textcolor{blue}{\textrm{on}}&\textcolor{blue}{\Gamma_{1}\times]0,+\infty[},
\\
\textcolor{blue}{u_{1}=0}&\textcolor{blue}{\textrm{on}}&\textcolor{blue}{\Gamma_{1}\times]0,+\infty[},
\\
\textcolor{red}{\Delta u_{2}=-a\,\p_{t}\pn u_{2}}&\textcolor{red}{\textrm{on}}&\textcolor{red}{\Gamma_{2}\times]0,+\infty[},
\\
\textcolor{red}{\pn\Delta u_{2}=b\,\p_{t} u_{2}}&\textcolor{red}{\textrm{on}}&\textcolor{red}{\Gamma_{2}\times]0,+\infty[},
\\
u_{1}(x,0)=u_{1}^{0}(x),\;\p_{t}u_{1}(x,0)=u_{1}^{1}(x)&\textrm{in}&\Omega_{1},
\\
u_{2}(x,0)=u_{2}^{0}(x),\;\p_{t}u_{2}(x,0)=u_{2}^{1}(x)&\textrm{in}&\Omega_{2}.
\end{array}\right.
\end{equation}
Where $\nu$ denotes the inner unit normal to the boundary, $c_{1}$, $c_{2}$ are constants strictly positives and $a$ and $b$ are a non negative bounded functions on $\Gamma_{2}$. Then we can suppose that there exists a strictly positive constant $c_{0}$ such that
\begin{equation}\label{2}
a\geq c_{0} \quad\textrm{and}\quad b\geq c_{0}\;\;\;\;\textrm{ on }\Gamma_{2}.
\end{equation}
And finally as regards $u_{1}^{0}$, $u_{1}^{1}$, $u_{2}^{0}$ and $u_{2}^{1}$ we will fix them later in the right spaces.

The energy of a solution
\begin{equation*}
u=\left\{
\begin{array}{lcl}
u_{1}&\textrm{in}&\Omega_{1}
\\
u_{2}&\textrm{in}&\Omega_{2}
\end{array}\right.
\end{equation*}
of~\eqref{1} at the time $t\geq 0$ is defined by
\begin{equation*}
E(t,u)=\frac{1}{2}\int_{\Omega}\Big(|\p_{t}u(x,t)|^{2}+\alpha^{2}(x)|\Delta u(x,t)|^{2}\Big)\alpha^{-1}(x)\,\ud x,
\end{equation*}
where
\begin{equation*}
\alpha=\left\{
\begin{array}{lcl}
c_{1}&\textrm{in}&\Omega_{1},
\\
c_{2}&\textrm{in}&\Omega_{2}.
\end{array}\right.
\end{equation*}
By Green's formula we can prove that for all $\;t_{1},\,t_{2}>0$ we have
\[
E(t_{2},u)-E(t_{1},u)=-c_{2}\int_{t_{1}}^{t_{2}}\!\!\!\int_{\Gamma_{2}}a|\p_{t}\pn u_{2}(x,t)|^{2}\,\ud x\,\ud t-c_{2}\int_{t_{1}}^{t_{2}}\!\!\!\int_{\Gamma_{2}}b|\p_{t}u_{2}(x,t)|^{2}\,\ud x\,\ud t,
\]
and this mean that the energy is decreasing over the time.

We define the operator $\mathcal{A}$ by
$$\mathcal{A}=\left[\begin{array}{cc}
0& \id
\\
-\alpha^{2}\Delta^{2}& 0
\end{array}\right]$$
the Hilbert space $\mathcal{H}=X\times H$ where $H=L^{2}(\Omega,\alpha^{-1}(x)\,\ud x)$ and
\begin{equation}\label{8}
\begin{split}
X=\big\{u\in H\,:\,u_{1}\in H^{2}(\Omega_{1}),\,u_{2}\in H^{2}(\Omega_{2}),\,u_{1\,|\Gamma_{1}}=0,\,u_{1\,|\Gamma_{0}}=u_{2\,|\Gamma_{0}},
\\
\pn u_{1\,|\Gamma_{0}}=\pn u_{2\,|\Gamma_{0}},\int_{\Gamma_{2}}u_{2}.\overline{\pn u}_{2}\ud x=0\big\},
\end{split}
\end{equation}
and the domain of $\mathcal{A}$ by
\begin{align*}
\mathcal{D}(\mathcal{A})=\big\{&(u,v)\in\mathcal{H}\,:\,(v,\Delta^{2}u)\in\mathcal{H},\,c_{1}\Delta u_{1\,|\Gamma_{0}}=c_{2}\Delta u_{2\,|\Gamma_{0}},\,\Delta u_{1\,|\Gamma_{1}}=0,
\\
&c_{1}\pn \Delta u_{1\,|\Gamma_{0}}=c_{2}\pn \Delta u_{2\,|\Gamma_{0}},\,\Delta u_{2\,|\Gamma_{2}}=-a\,\pn v_{2\,|\Gamma_{2}},\,\pn\Delta u_{2\,|\Gamma_{2}}=b\,v_{2\,|\Gamma_{2}}\big\}.
\end{align*}
Now we are able to state our main results
\begin{thm}\label{9}
There exists $C_{1},\,C_{2},\,C_{3}>0$ such that if $|\im(\lambda)|\leq C_{1}\e^{-C_{2}|\re(\lambda)|}$ and $|\lambda|>C_{3}$ the resolvent $(\lambda\id+i\mathcal{A})^{-1}$ is analytic and moreovere we have
$$\|(\lambda\id+i\mathcal{A})^{-1}\|_{\mathcal{L}(\mathcal{H},\mathcal{H})}\leq C\e^{C|\re(\lambda)|}.$$
\end{thm}
As an immediate consequence (see~\cite[p.17]{Bur} and also more recently~\cite{BD}) of the previous theorem, we get the following rate of decrease of energy
\begin{thm}\label{10}
For any $k>0$ there exists $C>0$ such that for any initial data $(u_{0},v_{0})=(u_{1}^{0},u_{2}^{0},u_{1}^{1},u_{2}^{1})\in\mathcal{D}(\mathcal{A}^{k})$ the solution $u(x,t)$ of~\eqref{1} starting from $(u_{0},v_{0})$ satisfy
$$E(t,u)\leq\frac{C}{(\ln(2+t))^{2k}}\|(u_{0},v_{0})\|_{\mathcal{D}(\mathcal{A}^{k})}^{2},\quad \forall\; t>0.$$
\end{thm}
\begin{rem}\rm{
\*
\begin{enumerate}
	\item[1)]In the case where $\Gamma_{1}=\emptyset$, (i.e $\Gamma=\Gamma_{2}$) and we have only one boundary term effective in all $\Gamma$ gived by $\Delta u_{2\,|\Gamma}=-a\,\p_{t}\pn u_{2\,|\Gamma}$, Ammari and Vodev~\cite{AV} have proved an exponential stabilization result to the system~\eqref{1}.
	\item[2)]Theorem~\ref{10} remains valid if we suppose that $\Gamma_{1}=\emptyset$ and $\Gamma_{2}=\Gamma_{2}'\cup\Gamma_{2}''$ satisfying $\Gamma_{2}'\cap\Gamma_{2}''=\emptyset$ with the following transmission and boundary conditions:
	\begin{equation*}
\left\{
\begin{array}{lll}
\textcolor{magenta}{u_{1}=u_{2}}&\textcolor{magenta}{\textrm{on}}&\textcolor{magenta}{\Gamma_{0}\times]0,+\infty[},
\\
\textcolor{magenta}{\pn u_{1}=\pn u_{2}}&\textcolor{magenta}{\textrm{on}}&\textcolor{magenta}{\Gamma_{0}\times]0,+\infty[},
\\
\textcolor{magenta}{c_{1}\Delta u_{1}=c_{2}\Delta u_{2}}&\textcolor{magenta}{\textrm{on}}&\textcolor{magenta}{\Gamma_{0}\times]0,+\infty[},
\\
\textcolor{magenta}{c_{1}\pn \Delta u_{1}=c_{2}\pn \Delta u_{2}}&\textcolor{magenta}{\textrm{on}}&\textcolor{magenta}{\Gamma_{0}\times]0,+\infty[},
\\
\textcolor{blue}{\Delta u_{1}=0}&\textcolor{blue}{\textrm{on}}&\textcolor{blue}{\Gamma_{2}'\times]0,+\infty[},
\\
\textcolor{blue}{u_{1}=0}&\textcolor{blue}{\textrm{on}}&\textcolor{blue}{\Gamma_{2}'\times]0,+\infty[},
\\
\textcolor{red}{\Delta u_{2}=-a\,\p_{t}\pn u_{2}}&\textcolor{red}{\textrm{on}}&\textcolor{red}{\Gamma_{2}''\times]0,+\infty[},
\\
\textcolor{red}{\pn\Delta u_{2}=b\,\p_{t} u_{2}}&\textcolor{red}{\textrm{on}}&\textcolor{red}{\Gamma_{2}''\times]0,+\infty[}.
\end{array}\right.
\end{equation*}
	\item[3)]To prove Theorem~\ref{9} and Theorem~\ref{10}, we make use the Carleman estimates to obtain information about the resolvent in a boundary domain, the cost is to use phases functions satisfying H\"ormander's assumption. Albano~\cite{A} proved a Carleman estimate for the plate operator, by decomposing the operator as the product of two Schr{\"o}dinger ones and gives for eatch of them the corresponding Carleman estimate then by making together these two estimates we obtain the result. Inspired from this method, we are going to do one similar decomposition of a system of forth order to get two new systems of second order, then we apply the classical Carleman estimate to one of these derivative operators to obtain easly the result of Theorem~\ref{9}.
	\item[4)]Theorem~\ref{9} and Theorem~\ref{10} are analogous to those of Fathallah~\cite{I}, in the case of hyperbolic-parabolic coupled system, and Lebeau and Robbiano~\cite{LR} resuts, in the case of scalar wave equation without transmission, but our method is different from their because it consist to use the Carleman estimates directly for the stationary operator without going through the interpolation inequality.
	\item [5)]Several studies have focused on transmission problems, such as the works of Bellassoued~\cite{B} and Fathallah~\cite{I} for the stabilization problems and that of Le Rousseau and Robbiano~\cite{RR} for a control problem. In eatch of these works we need to find a Carleman estimates near the interface, but here and thanks to the transmission conditions we will use only the classical Carleman estimates (See for instance Le Rouseau and Lebeau in~\cite{LL} and Lebeau and Robbiano in~\cite{LR} and~\cite{LR2}).
\end{enumerate}
}
\end{rem}

In this paper $C$ will always be a generic positive constant whose value may be different from one line to another.

The outline of this paper is as follow. In section~\ref{11} we prove the well-Posedness of the problem~\eqref{1} and in section~\ref{d1} we prove the logarithmic decay of energy of the system~\eqref{1}.
\section{Well-Posedness of the problem}\label{11}
To prove the Well-Posedness of the problem~\eqref{1} we are going to use the semigroups theory. Our strategy consiste to write the equations as a Cauchy problem with an operator which generates a semigroup of contractions.
\subsection{The Cauchy problem}
Throughout this paper, we denote $\mathcal{O}=\Omega_{1}\cup\Omega_{2}$ and $\langle\,.\,,\,.\,\rangle_{H}$ the inner product in $H=L^{2}(\mathcal{O},\alpha^{-1}(x)\,\ud x)$ defined by
\begin{equation*}
\langle u,v\rangle_{H}=\int_{\mathcal{O}}u(x)\overline{v(x)}\alpha^{-1}(x)\,\ud x=\int_{\Omega_{1}}u_{1}(x)\overline{v_{1}(x)}c_{1}^{-1}\,\ud x+\int_{\Omega_{2}}u_{2}(x)\overline{v_{2}(x)}c_{2}^{-1}\,\ud x,
\end{equation*}
where we recall that
$$u(x,t)=\left\{\begin{array}{lcl}
u_{1}(x,t)&\textrm{if}&x\in\Omega_{1},
\\
u_{2}(x,t)&\textrm{if}&x\in\Omega_{2},
\end{array}\right.
\quad \text{and} \quad
v(x,t)=\left\{\begin{array}{lcl}
v_{1}(x,t)&\textrm{if}&x\in\Omega_{1},
\\
v_{2}(x,t)&\textrm{if}&x\in\Omega_{2}.
\end{array}\right.$$
The tow first equations of~\eqref{1} can be written as follows
$$\p_{t}\left(\begin{array}{c}
u
\\
v
\end{array}\right)=\mathcal{A}\left(\begin{array}{c}
u
\\
v
\end{array}\right)$$
and the Cauchy problem is written in following form
$$\left\{\begin{array}{lcl}
\p_{t}\left(\begin{array}{c}
u
\\
v
\end{array}\right)(x,t)=\mathcal{A}\left(\begin{array}{c}
u
\\
v
\end{array}\right)(x,t)&\text{if}&(x,t)\in\mathcal{O}\times]0,+\infty[,
\\
\left(\begin{array}{l}
u
\\
v
\end{array}\right)(x,0)=\left(\begin{array}{l}
u_{0}
\\
v_{0}
\end{array}\right)(x)&\text{if}&x\in\mathcal{O},
\end{array}\right.$$
where
$$u_{0}(x)=\left\{\begin{array}{lcl}
u_{1}^{0}(x)&\textrm{if}&x\in\Omega_{1},
\\
u_{2}^{0}(x)&\textrm{if}&x\in\Omega_{2},
\end{array}\right.
\quad \text{and} \qquad
v_{0}(x)=\left\{\begin{array}{lcl}
v_{1}^{0}(x)&\textrm{if}&x\in\Omega_{1},
\\
v_{2}^{0}(x)&\textrm{if}&x\in\Omega_{2}.
\end{array}\right.$$

Now we have to specify the functional space and the domain of the operator $\mathcal{A}$, that's why we are going first to define anothor operator $G$ which is the square root of the operator $\alpha^{2}\Delta^{2}$ (see~\cite[p.391]{TucWei} for the definition and more details).
\subsection{Properties of the square root of the operator $\alpha^{2}\Delta^{2}$}
In the space $H=L^{2}(\mathcal{O},\alpha^{-1}(x)\,\ud x)$ we define the operator $G$ by the following expression
\[
Gu=-\alpha\Delta u\qquad \forall\;u\in\mathcal{D}(G)
\]
with domain $\mathcal{D}(G)=X$ defined in~\eqref{8}. The space $X$ is equipped with the norm $\|u\|_{X}=\|Gu\|_{H}$ and we defined the graph norm of $G$ by
$$\|u\|_{gr(G)}^{2}=\|u\|_{H}^{2}+\|Gu\|_{H}^{2}$$
then we have the following result
\begin{pro}
$(X,\|\,.\,\|_{X})$ is a Hilbert space with a norm equivalent to the graph norm of $G$.
\end{pro}
\begin{pr}
It is easy to show that if $G$ is a colsed operator then $(X,\|\,.\,\|_{gr(G)})$ is a Hilbert space. Thus to prove the proposition it suffices to show that $G$ is closed and both norms are equivalent.
\\
By Green's formula and Poincar\'e inequality it is easy to show that there exists $C>0$ such that
$$\langle Gu,u\rangle_{H}=\|\nabla u\|_{L^{2}(\mathcal{O})}^{2}\geq C\|u\|_{H}^{2}\quad\forall\,u\in X$$
then $G$ is a strictly positive operator and we have
$$\|Gu\|_{H}\|u\|_{H}\geq\langle Gu,u\rangle_{H}\geq C\|u\|_{H}^{2}\quad\forall\,u\in X$$
which prove the equivalence between the tow norms.
\\
Since $G$ is positive then by Propsition 3.3.5 in~\cite[p.79]{TucWei} $-G$ is m-dissipative and thus $G$ is a closed operator. This completes the proof.
\end{pr}

This last result allows us to properly define the functional space of the operator $\mathcal{A}$.
\begin{pro}\label{7}
The two spaces $(X,\|\,.\,\|_{2})$ and $(X,\|\,.\,\|_{X})$ are algebraically and topologically the same. Where $\|\,.\,\|_{2}$ is the classical Sobolev norm.
\end{pro}
\begin{pr}
We have only to prove that the two norms are equivalent.
\\
First, we note that $(X,\|\,.\,\|_{2})$ is a Hilbert space because $X$ is a closed subspace of $H^{2}(\mathcal{O})$, in addition we have
$$\|u\|_{X}=\|Gu\|_{L^{2}(\mathcal{O})}\leq C\|u\|_{2}\quad\forall\,u\in X,$$
and while $(X,\|\,.\,\|_{X})$ is also a Hilbert space, then according to the Banach theorem (see Corollary 9.2.3 from~\cite[p.132]{YVA}) the tow norms are equivalent.
\end{pr}

As an important consequence of this result is that the space $X$ is a Hilbert space with the norm $\|\,.\,\|_{gr(G)}+\|a^{\frac{1}{2}}\pn\,.\,\|_{L^{2}(\Gamma_{2})}+\|b^{\frac{1}{2}}\,.\,\|_{L^{2}(\Gamma_{2})}$.
\subsection{Existence and uniqueness of the solution}
We set $\mathcal{H}=X\times H$ the Hilbert space with the norm
\[
\|(u,v)\|=\|u\|_{X}+\|v\|_{H}\qquad\forall\, (u,v)\in\mathcal{H},
\]
and we recall that the domain of the operator $\mathcal{A}$ is defined by
\begin{align*}
\mathcal{D}(\mathcal{A})=\big\{&(u,v)\in\mathcal{H}\,:\,(v,\Delta^{2}u)\in\mathcal{H},\,c_{1}\Delta u_{1\,|\Gamma_{0}}=c_{2}\Delta u_{2\,|\Gamma_{0}},\,\Delta u_{1\,|\Gamma_{1}}=0,
\\
&c_{1}\pn \Delta u_{1\,|\Gamma_{0}}=c_{2}\pn \Delta u_{2\,|\Gamma_{0}},\,\Delta u_{2\,|\Gamma_{2}}=-a\,\pn v_{2\,|\Gamma_{2}},\,\pn\Delta u_{2\,|\Gamma_{2}}=b\,v_{2\,|\Gamma_{2}}\big\}.
\end{align*}
\begin{thm}
Under the above assumptions, $\mathcal{A}$ is m-dissipative and especially it generates a strongly semigroup of contractions in $\mathcal{H}$.
\end{thm}
\begin{pr}
According to Lumer-Phillips theorem (see for exemple~\cite[p.103]{TucWei}) we have only to prove that $\mathcal{A}$ is m-dissipative.
\\
Let $(u,v)\in\mathcal{D}(\mathcal{A})$ then by Green's formula we have
\begin{eqnarray*}
\re\left(\left\langle\mathcal{A}\left(\begin{array}{l}
u
\\
v
\end{array}\right),\left(\begin{array}{l}
u
\\
v
\end{array}\right)\right\rangle_{\mathcal{H}}\right)\!\!\!\!&=&\!\!\!\!\alpha\,\re\left(\langle\Delta u,\Delta v\rangle_{L^{2}(\mathcal{O})}-\langle\Delta^{2}u,v\rangle_{L^{2}(\mathcal{O})}\right)
\\
\!\!\!\!&=&\!\!\!\!\alpha\,\re\left(\langle\Delta u,\pn v\rangle_{L^{2}(\p\mathcal{O})}-\langle\pn \Delta u,v\rangle_{L^{2}(\p\mathcal{O})}\right)
\\
\!\!\!\!&=&\!\!\!\!-c_{2}\|a^{\frac{1}{2}}\pn v_{2}\|_{L^{2}(\Gamma_{2})}^{2}-c_{2}\|b^{\frac{1}{2}} v_{2}\|_{L^{2}(\Gamma_{2})}^{2}\leq 0.
\end{eqnarray*}
This shows that $\mathcal{A}$ is dissipative.
\\
Let now $(f,g)\in\mathcal{H}$ and our purpose is to find a couple $(u,v)\in\mathcal{D}(\mathcal{A})$ such that
$$\left(\id-\mathcal{A}\right)\left(\begin{array}{l}
u
\\
v
\end{array}\right)=\left(\begin{array}{l}
u-v
\\
v+\alpha^{2}\Delta^{2}u
\end{array}\right)=\left(\begin{array}{l}
f
\\
g
\end{array}\right)$$
and more explicitly we have to find $(u,v)\in\mathcal{D}(\mathcal{A})$ such that
\begin{equation*}
\left\{\begin{array}{l}
v=u-f=\left\{\begin{array}{lll}
v_{1}=u_{1}-f_{1}&\text{in}&\Omega_{1}
\\
v_{2}=u_{2}-f_{2}&\text{in}&\Omega_{2}
\end{array}\right.
\\
u+\alpha^{2}\Delta^{2} u=f+g=\left\{\begin{array}{lll}
u_{1}+c_{1}^{2}\Delta^{2} u_{1}=f_{1}+g_{1}&\text{in}&\Omega_{1}
\\
u_{2}+c_{2}^{2}\Delta^{2} u_{2}=f_{2}+g_{2}&\text{in}&\Omega_{2}.
\end{array}\right.
\end{array}\right.
\end{equation*}
First note that, thanks to the remark after Proposition~\ref{7} and the Riesz representation theorem, there exists a unique $u\in X=\mathcal{D}(G)$ such that for all $\varphi\in X$ we have
\begin{equation}\label{4}
\begin{split}
\langle f+g,\varphi\rangle_{H}+\langle c_{2}a\,\pn f_{2},\pn\varphi_{2}\rangle_{L^{2}(\Gamma_{2})}+\langle c_{2}b\,f_{2},\varphi_{2}\rangle_{L^{2}(\Gamma_{2})}=\langle\alpha\Delta u,\alpha\Delta\varphi\rangle_{H}
\\
+\langle u,\varphi\rangle_{H}+\langle c_{2}a\,\pn u_{2},\pn\varphi_{2}\rangle_{L^{2}(\Gamma_{2})}+\langle c_{2}b\, u_{2},\varphi_{2}\rangle_{L^{2}(\Gamma_{2})}.
\end{split}
\end{equation}
In particular for all $\varphi\in\Cc(\mathcal{O})$ the expression~\eqref{4} is written as follows
$$\langle\alpha\Delta^{2}u+\alpha^{-1}(u-f-g),\varphi\rangle_{L^{2}(\mathcal{O})}=0$$
then we have
\begin{equation}\label{5}
u+\alpha^{2}\Delta^{2} u=f+g\quad\text{in}\;L^{2}(\mathcal{O}).
\end{equation}
Now if we return again to the expression~\eqref{4} then through Green's formula we write it as follows
\begin{align*}
&\langle\alpha\Delta^{2}u+\alpha^{-1}(u-f-g),\varphi\rangle_{L^{2}(\mathcal{O})}=\langle c_{2}a\,\pn f_{2},\pn\varphi_{2}\rangle_{L^{2}(\Gamma_{2})}+\langle c_{2}b\,f_{2},\varphi_{2}\rangle_{L^{2}(\Gamma_{2})}+
\\
&\langle\alpha\pn\Delta u,\varphi\rangle_{L^{2}(\p\mathcal{O})}-\langle c_{2}a\,\pn u_{2},\pn\varphi_{2}\rangle_{L^{2}(\Gamma_{2})}-\langle c_{2}b\,u_{2},\varphi_{2}\rangle_{L^{2}(\Gamma_{2})}-\langle\alpha\Delta u,\pn\varphi\rangle_{L^{2}(\p\mathcal{O})}
\end{align*}
then by~\eqref{5} and after a simple calculation we get that for all $\varphi\in X$ that
\begin{align*}
&\langle c_{1}\pn\Delta u_{1}-c_{2}\pn\Delta u_{2},\varphi_{2}\rangle_{L^{2}(\Gamma_{0})}-\langle c_{1}\Delta u_{1}-c_{2}\Delta u_{2},\pn\varphi_{1}\rangle_{L^{2}(\Gamma_{0})}-\langle c_{1}\Delta u_{1},\pn\varphi_{1}\rangle_{L^{2}(\Gamma_{1})}
\\
&+c_{2}\langle a(\pn f_{2}-\pn u_{2})-\Delta u_{2},\pn\varphi_{2}\rangle_{L^{2}(\Gamma_{2})}+c_{2}\langle b(f_{2}-u_{2})+\pn\Delta u_{2},\varphi_{2}\rangle_{L^{2}(\Gamma_{2})}=0
\end{align*}
and this shows the following equalities
$$c_{1}\Delta u_{1\,|\Gamma_{0}}=c_{2}\Delta u_{2\,|\Gamma_{0}},\; c_{1}\pn \Delta u_{1\,|\Gamma_{0}}=c_{2}\pn \Delta u_{2\,|\Gamma_{0}},\; \Delta u_{1\,|\Gamma_{1}}=0$$
and also
\begin{equation*}
\begin{array}{l}
\Delta u_{2\,|\Gamma_ {2}}=-a(\pn u_{2\,|\Gamma_ {2}}-\pn f_{2\,|\Gamma_ {2}})=-a\,\pn v_{2\,|\Gamma_ {2}}\quad
\\
\pn\Delta u_{2\,|\Gamma_ {2}}=b(u_{2\,|\Gamma_ {2}}- f_{2\,|\Gamma_ {2}})=b\,v_{2\,|\Gamma_ {2}}.
\end{array}
\end{equation*}
And this concludes the proof.
\end{pr}

One consequence of this last result is that if we assume that $(u_{0},v_{0})\in\mathcal{D}(\mathcal{A})$, there exists a unique solution of~\eqref{1} which can be expressed by means of a semigroup on $\mathcal{H}$ as follows
\begin{equation}\label{6}
\left(\begin{array}{l}
u
\\
\p_{t}u
\end{array}\right)=e^{t\mathcal{A}}\left(\begin{array}{l}
u_{0}
\\
v_{0}
\end{array}\right)
\end{equation}
where $e^{t\mathcal{A}}$ is the semigroupe of the operator $\mathcal{A}$. And we have the following regularity of the solution
$$
\left(\begin{array}{l}
u
\\
\p_{t}u
\end{array}\right)
\in C([0,+\infty[,\mathcal{D}(\mathcal{A}))\cap C^{1}([0,+\infty[,\mathcal{H}).$$
\\
And if $(u_{0},v_{0})\in\mathcal{H}$, the function $u(t)$ given by~\eqref{6} is the mild solution of~\eqref{1} and it lives in $C([0,+\infty[,\mathcal{H})$.
\section{Proof of Theorem~\ref{9}}\label{d1}
The purpose of this section is to find an estimate of the resolvent $(\lambda\id+i\mathcal{A})^{-1}$ for $\lambda$ in the region $\{z\in\mathbb{C};\;|\im(z)|<C_{1}\e^{-C_{2}|\re(z)|},\,|z|>C_{3}\}$ with some constants $C_{1},\,C_{2},\,C_{3}>0$. More precisely we prove that $\|(\lambda\id+i\mathcal{A})^{-1}\|_{\mathscr{L}(\mathcal{H},\mathcal{H})}\leq C\e^{C|\re(\lambda)|}$ which imply the weak energy decay of the solution of the equation~\eqref{1}.

The main idea consiste to the use of the Carleman estimates for a second order elliptic operator which it derived from an original one of fourth order and this is what comes from the originality of our work, it means we prove the stability result for a system of fourth order by using an estimate of Carleman of second order only.

As manshed previously that to prove Theorem~\ref{9}, we will need the Carleman estimates due to Lebeau and Robbiano~\cite{LR2} and formuled by Burq~\cite{Bur}. We consider the elliptic second order operator $P=-h^{2}\Delta$ defined for a complex valued functions which are defined in an open subset $U\subset\Rn$ with smooth boundary, and whose principal symbol is denoted by $p(x,\xi)=|\xi|^{2}$, where $h$ is a very small semi-classical parameter. 

Let $\vfi\in\Ci(\overline{U})$ a real value function and let's define the adjoint operator $P_{\vfi}=\e^{\vfi/h}P\e^{-\vfi/h}$ of principal symbol $p_{\vfi}(x,\xi)=p(x,\xi+i\nabla\vfi)$ for $0<h\leq h_{0}$. Then we have the following result
\begin{pro}\cite[Proposition 2]{LR2}~\cite[Proposition 1]{LR}\label{d16}
Let $\gamma$ be an non-empty union of connex component of $\p U$. Assume the weight function $\vfi$ satisfies to the following assumptions:
\begin{enumerate}
	\item $\nabla\vfi\neq 0$ for all $x\in\overline{U}$
	\item $\pn\vfi\neq 0$ for all $x\in\p U$
	\item $\pn\vfi<0$ for all $x\in \gamma$
	\item The H{\"o}rmander's sub-ellipticity condition
	\[
	\forall\; (x,\xi)\in\overline{U}\times\Rn;\; p_{\varphi}(x,\xi)=0\Longrightarrow\{\re(p_{\varphi}),\im(p_{\varphi})\}(x,\xi)>0.
	\]
\end{enumerate}
Then there exists $C>0$ such that for all $u\in\Ci(\overline{U})$ satisfying
\[
\left\{\begin{array}{ll}
\displaystyle\Delta u=f &\text{in }U
\\
u=0&\text{on }\gamma,
\end{array}
\right.
\]
and for all $h\in]0,h_{0}]$ small we have
\begin{equation}\label{d6}
\begin{split}
h\int_{U}\e^{2\vfi/h}|u|^{2}\,\ud x&+h^{3}\int_{U}\e^{2\vfi/h}|\nabla u|^{2}\,\ud x\leq C\Big(h^{4}\int_{U}\e^{2\vfi/h}|f|^{2}\,\ud x
\\
&+h\int_{\p U\backslash\gamma}\e^{2\vfi/h}|u|^{2}\,\ud x+h^{3}\int_{\p U\backslash\gamma}\e^{2\vfi/h}|\pn u|^{2}\,\ud x\Big).
\end{split}
\end{equation}
\end{pro}
\begin{rem}
\rm{
\*
\begin{enumerate}
	\item [1)]If  the function $u$ is supported away from a subset $\gamma_{0}\subset\p U$ then the estimate~\eqref{d6} is allows true even if we don't assume that $\pn\vfi\neq 0$ in $\gamma_{0}$, while the proof is local.
	\item [2)]We can not assume that $\pn\vfi<0$ in the whole $\p U$, otherwise the weight function attain his global maximum in $U$, and thus our srtategy of the construction of the phases is fails (See next subsection).
\end{enumerate}
}
\end{rem}
\subsection{Weight function's construction}\label{d15}
\begin{figure}[htbp]
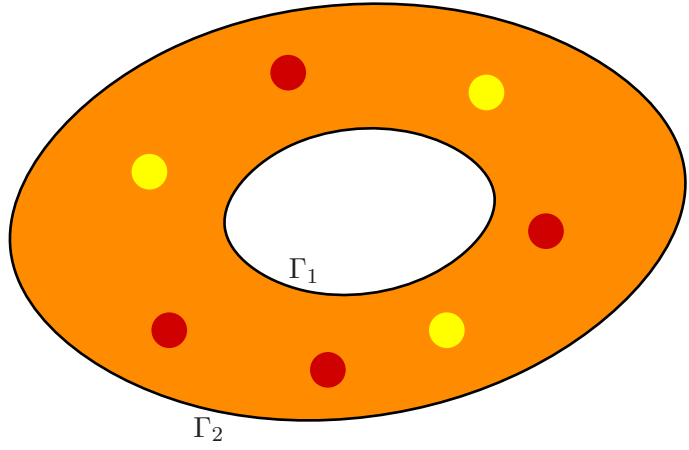

\figinit{1.5pt}
\figpt 0:(20,0)
\figpt 1:(-70,-10)
\figpt 2:(-20,-50)
\figpt 3:(60,-40)
\figpt 4:(100,10)
\figptsym 5:=3/1,4/
\figptsym 6:=2/1,4/
\figpthom 14:=0/1,0.6/
\figpthom 15:=0/2,0.6/
\figpthom 16:=0/3,0.6/
\figpthom 17:=0/4,0.6/
\figpthom 18:=0/5,0.6/
\figpthom 19:=0/6,0.6/
\figpt 7:(50,30)
\figpt 8:(0,35)
\figpt 9:(-35,10)
\figpt 10:(-30,-30)
\figpt 11:(65,-5)
\figpt 12:(10,-40)
\figpt 13:(40,-30)
\psbeginfig{}
\psset(width=1)
\psset(fillmode=yes,color=\DarkOrangergb)
\pscurve[1,2,3,4,5,6,1,2,3]
\psset(fillmode=yes,color=0.6)
\psset(fillmode=yes,color=1)
\pscurve[14,15,16,17,18,19,14,15,16]
\psset(color=1)
\pscirc 7(4)
\pscirc 8(4)
\pscirc 9(4)
\pscirc 10(4)
\pscirc 11(4)
\pscirc 12(4)
\pscirc 13(4)
\psset(fillmode=no)
\psset(color=\Blackrgb)
\pscurve[1,2,3,4,5,6,1,2,3]
\pscurve [14,15,16,17,18,19,14,15,16]
\psset(fillmode=yes,color=0.8\Redrgb)
\pscirc 8(4.5)
\pscirc 10(4.5)
\pscirc 11(4.5)
\pscirc 12(4.5)
\psset(color=\Yellowrgb)
\pscirc 7(4.5)
\pscirc 9(4.5)
\pscirc 13(4.5)
\psendfig
\figvisu{\figBoxA}{}{
\figwriten 15:$\Gamma_{1}$(2)
\figwrites 2:$\Gamma_{2}$(2)
}
\centerline{\box\figBoxA}
\centerline{\box\figBoxA}
\caption{The domains of the weight functions $\varphi_{1}$ and $\psi_{1}$ (in yellow and orange), $\varphi_{2}$ and $\psi_{2}$ (in red and orange) where they have not critical points.}
\label{fig2}
\end{figure}

In this section we will try to find two phases $\vfi_{1}$ and $\vfi_{2}$ which satisfy to the H{\"o}rmander's condition except in a finite number of ball where one of them do not satisfies this condition the second does and is strictly greater. The main ingredient of this section is the following one. Note that this result is similar to the Burq's one~\cite[Proposition 3.2]{Bur}, but here we give a new proof due to F.~Laudenbach.
\begin{pro}\label{d8}
With keeping the same notations as the first section, then there exists two real functions $\psi_{1},\,\psi_{2}\in\Ci(\Omega)$ satisfying for $k=1,2$ that $\pn\psi_{k\,|\Gamma}\neq 0$ and $\pn\psi_{k\,|\Gamma_{1}}<0$ having only degenerate critical points (of finite number) such that when $\nabla\psi_{k}=0$ then $\nabla\psi_{\sigma(k)}\neq 0$ and $\psi_{\sigma(k)}>\psi_{k}$. Where $\sigma$ is the permutation of the set $\{1,2\}$ different from the identity. 
\end{pro}
\begin{rem}
\rm{
\*
\begin{enumerate}
	\item[1)] One consequence of Proposition~\ref{d8} is that there exists a finite number of points $x_{kj_{k}}$ for $k=1,2$ and $j_{k}=1,\ldots,N_{k}$ and $\epsilon>0$ such that $B(x_{kj_{k}},2\epsilon)\subset\overline{\Omega}$ and $B(x_{1j_{1}},2\epsilon)\cap B(x_{2j_{2}},2\epsilon)=\emptyset$, for all $k=1,2$ and $j_{k}=1,\ldots,N_{k}$ and in $B(x_{kj_{k}},2\epsilon)$ we have $\psi_{\sigma(k)}>\psi_{k}$ (See Figure~\ref{fig2}).
	\item[2)] For $\lambda>0$ large enough the weight functions $\vfi_{k}=\e^{\lambda\psi_{k}}$ satisfy the H{\"o}rmander's condition in $\displaystyle U_{k}=\Omega\bigcap\left(\bigcup_{j_{k}=1}^{N_{k}}B(x_{kj_{k}},\epsilon)\right)^{c}$. Indeed, we have only to prove that for an open bounded subset $U\in\Rn$ and if $\psi\in\Ci(\overline{U})$ satisfying $|\nabla\psi|\geq C$ in $\overline{U}$ and $\vfi=\e^{\lambda\psi}$ we have $\{\re(p_{\vfi}),\im(p_{\vfi})\}(x,\xi)\geq C'$ in $\overline{U}\times\Rn$ for $\lambda>0$ large enough. We have
$$
\left\{\begin{array}{c}
\nabla\vfi=\lambda\e^{\lambda\psi}\nabla\psi\;\text{ and }\;\vfi''=\e^{\lambda\psi}(\lambda\nabla\psi.{}^{t}\nabla\psi+\lambda\psi'')
\\
p_{\vfi}(x,\xi)=0\Longrightarrow\langle\xi,\nabla\vfi\rangle=0\text{ and }|\xi|^{2}=|\nabla\vfi|^{2}
\end{array}\right.
$$
then we obtain
\begin{eqnarray*}
\{\re(p_{\vfi}),\im(p_{\vfi})\}(x,\xi)&=&4\lambda\e^{\lambda\psi}\,{}^{t}\xi.\psi''.\xi+4\e^{3\lambda\psi}(\lambda^{4}|\nabla\psi|^{2}+\lambda^{3}\,{}^{t}\nabla\psi.\psi''.\nabla\psi)
\\
&=&4\e^{3\lambda\psi}(\lambda^{4}|\nabla\psi|^{2}+O(\lambda^{3})).
\end{eqnarray*}
Which conclude the result.
	\item[3)] In general, Proposition~\ref{d8} is also true for any smooth manifold with boundary which the latter is the disjoint union of two open and closed submanifolds.
\end{enumerate}
}
\end{rem}
\begin{pr}
While the Morse functions are dense (for the $\Ci$ topology) in the set of $\Ci$ functions then we can find $\psi_{1}$ a Morse function such that $\pn\psi_{1\,|\Gamma_{1}}<0$ and $\pn\psi_{1\,|\Gamma_{2}}>0$. We can suppose that $\psi_{1}$ have no local maximum in $\Omega$ (The proceeding of the elimination of the maximum is described by Burq~\cite[Appendix A]{Bur}, we can see also~\cite[Theorem 8.1]{Mi} and~\cite[Lemma 2.6]{L}).

Let $c$ be a critical point of $\psi_{1}$ while its index is different from $n$ then we can find a $\Ci$ arc $\gamma_{c}:[-1,1]\rightarrow\Omega$ such that $\gamma_{c}(0)=c$ and $\psi_{1}(\gamma_{c}(1))=\psi_{1}(\gamma_{c}(-1))>\psi_{1}(c)$. We do this construction for all the critical points of $\psi_{1}$ so that all the arcs are mutually disjoint. Hence, this allows us to find a vector field $X$ in $\Omega$, vanishing near the boundary of $\Omega$ such that for all critical points $c$ of $\psi_{1}$ we have
$$X(\gamma_{c}(t))=\stackrel{.}{\gamma}_{c}(t),$$
where $\stackrel{.}{\gamma}$ stand for the time derivative.

We denote $\phi_{t}$ its flow:
$$\stackrel{.}{\phi_{t}}(x)=X(\phi_{t}(x)),$$
and we set $\psi_{2}=\psi_{1}\circ\phi_{1}$, thus $\psi_{1}$ and $\psi_{2}$ satisfy the required properties. Indeed, since $X\equiv0$ near the boundary $\Gamma$ which mean that $\phi_{t}(x)=x$ near $\Gamma$ then $\pn\psi_{1\,|\Gamma}=\pn\psi_{2\,|\Gamma}$. If $c$ is a critical point of $\psi_{1}$ then we have $\psi_{2}(c)=\psi_{1}(\gamma_{c}(1))>\psi_{1}(c)$, and if $c'$ is a critical point of $\psi_{2}$ then $c'=\phi_{-1}(c)$ where $c$ is a critical point of $\psi_{1}$ and we have $\psi_{2}(c')=\psi_{1}(\phi_{1}\circ\phi_{-1}(c))=\psi_{1}(c)<\psi_{1}(\phi_{-1}(c))=\psi_{1}(c')$ by the construction of $\gamma_{c}$.
\end{pr}
\subsection{Back to the proof of Theorem~\ref{9}}
We return now to the main proof. Let $(F,G)\in\mathcal{H}$ with $F=(f_{1},f_{2})\in X$ and $G=(g_{1},g_{2})\in H$ and $(u,v)\in D(\mathcal{A})$ with $u=(u_{1},u_{2})$ and $v=(v_{1},v_{2})$ such that
$$(\lambda\id+i\mathcal{A})\left(\begin{array}{c}
u
\\
v
\end{array}\right)=\left(\begin{array}{c}
F
\\
G
\end{array}\right),$$
then we get the following boundary value problem
\begin{equation}\label{d2}
\left\{\begin{array}{ll}
\lambda u+iv=F&\text{in }\mathcal{O}
\\
-i\alpha^{2}\Delta^{2}u+\lambda v=G&\text{in }\mathcal{O}
\\
u_{1}=u_{2},\quad\pn u_{1}=\pn u_{2}&\text{on }\Gamma_{0}
\\
c_{1}\Delta u_{1}=c_{2}\Delta u_{2},\quad c_{1}\pn \Delta u_{1}=c_{2}\pn\Delta u_{2}&\text{on }\Gamma_{0}
\\
u_{1}=0&\text{on }\Gamma_{1}
\\
\Delta u_{1}=0&\text{on }\Gamma_{1}
\\
\Delta u_{2}=-a\,\pn v_{2}&\text{on }\Gamma_{2}
\\
\pn\Delta u_{2}=b\, v_{2}&\text{on }\Gamma_{2}.
\end{array}\right.
\end{equation}
Then the solution $(u,v)$ of~\eqref{d2} satisfies
\begin{equation}\label{d4}
\left\{\begin{array}{ll}
v=i\lambda u-iF&\text{in }\mathcal{O}
\\
(-\lambda^{2}+\alpha^{2}\Delta^{2})u=iG-\lambda F=\Phi&\text{in }\mathcal{O}
\\
u_{1}=u_{2},\quad\pn u_{1}=\pn u_{2}&\text{on }\Gamma_{0}
\\
c_{1}\Delta u_{1}=c_{2}\Delta u_{2},\quad c_{1}\pn \Delta u_{1}=c_{2}\pn\Delta u_{2}&\text{on }\Gamma_{0}
\\
u_{1}=0&\text{on }\Gamma_{1}
\\
\Delta u_{1}=0&\text{on }\Gamma_{1}
\\
\Delta u_{2}+i\lambda a\,\pn u_{2}=ia\,\pn f_{2}=\phi_{a}&\text{on }\Gamma_{2}
\\
\pn\Delta u_{2}-i\lambda b\,u_{2}=-ib\,f_{2}=\phi_{b}&\text{on }\Gamma_{2}.
\end{array}\right.
\end{equation}
Integrating by part we obtain
\begin{eqnarray}\label{d3}
\langle\Phi,u\rangle_{H}&=&\!\!\!\alpha^{-1}\int(-\lambda^{2}+\alpha^{2}\Delta^{2})u.\overline{u}\,\ud x=-\lambda^{2}\|u\|_{H}^{2}+\alpha\int\Delta^{2}u.\overline{u}\,\ud x\nonumber
\\
&=&\!\!\!-\lambda^{2}\|u\|_{H}^{2}+\alpha^{2}\|\Delta u\|_{H}^{2}-c_{2}\langle \Delta u_{2},\pn u_{2}\rangle_{L^{2}(\Gamma_{2})}+c_{2}\langle \pn\Delta u_{2},u_{2}\rangle_{L^{2}(\Gamma_{2})}\nonumber
\\
&=&\!\!\!-\lambda^{2}\|u\|_{H}^{2}+\alpha^{2}\|\Delta u\|_{H}^{2}-c_{2}\langle\phi_{a},\pn u_{2}\rangle_{L^{2}(\Gamma_{2})}+c_{2}\langle\phi_{b},u_{2}\rangle_{L^{2}(\Gamma_{2})}
\\
&+&\!\!\!i\lambda c_{2}\langle a\,\pn u_{2},\pn u_{2}\rangle_{L^{2}(\Gamma_{2})}+i\lambda c_{2}\langle\nonumber b\,u_{2},u_{2}\rangle_{L^{2}(\Gamma_{2})}.
\end{eqnarray}
Keeping only the imaginary part of~\eqref{d3} then we get
\begin{equation}\label{d9}
|\re(\lambda)|\int_{\Gamma_{2}} a|\pn u_{2}|^{2}+ b|u_{2}|^{2}\,\ud x\leq C(\|\Phi\|_{H}\|u\|_{H}+2|\re(\lambda)\im(\lambda)|^{2}\,\|u\|_{H}^{2}+\|f_{2}\|_{2}\|u_{2}\|_{2}).
\end{equation}
Now we return to the system~\eqref{d4} which can be recast as follows
\begin{equation}\label{d7}
\left\{\begin{array}{ll}
v=i\lambda u-iF&\text{in }\mathcal{O}
\\
\displaystyle\Big(\Delta-\mathrm{sig}\big(\re(\lambda)\big)\frac{\lambda}{\alpha}\Big)u=\alpha^{-1}w&\text{in }\mathcal{O}
\\
u_{1}=u_{2},\quad\pn u_{1}=\pn u_{2}&\text{on }\Gamma_{0}
\\
u_{1}=0&\text{on }\Gamma_{1}
\end{array}\right.
\end{equation}
and
\begin{equation}\label{d5}
\left\{\begin{array}{ll}
\displaystyle\Big(\Delta+\mathrm{sig}\big(\re(\lambda)\big)\frac{\lambda}{\alpha}\Big)w=\alpha^{-1}\Phi&\text{in }\mathcal{O}
\\
w_{1}=w_{2},\quad\pn w_{1}=\pn w_{2}&\text{on }\Gamma_{0}
\\
w_{1}=0&\text{on }\Gamma_{1}
\\
\displaystyle w_{2}=c_{2}\phi_{a}-ic_{2}\lambda a\,\pn u_{2}-\mathrm{sig}\big(\re(\lambda)\big)\lambda u_{2}&\text{on }\Gamma_{2}
\\
\displaystyle\pn w_{2}=c_{2}\phi_{b}+ic_{2}\lambda b\,u_{2}-\mathrm{sig}\big(\re(\lambda)\big)\lambda\pn u_{2}&\text{on }\Gamma_{2},
\end{array}\right.
\end{equation}
where
\begin{equation*}
w=\left\{
\begin{array}{lcl}
w_{1}&\textrm{in}&\Omega_{1}
\\
w_{2}&\textrm{in}&\Omega_{2}.
\end{array}\right.
\end{equation*}
To prove the resolvent estimate, we need the following result which is a consequence of the Carleman estimates introduced in the beginning of this section.
\begin{lem}
There exists $C>0$ such that for any $u$ and $w$ solution of~\eqref{d7} and~\eqref{d5} the following estimate holds:
\begin{equation}\label{d10}
\|u\|_{X}^{2}\leq C\e^{C/h}\Big(\|\Phi\|_{L^{2}(\Omega)}^{2}+\|\phi_{a}\|_{L^{2}(\Gamma_{2})}^{2}+\|\phi_{b}\|_{L^{2}(\Gamma_{2})}^{2}+\int_{\Gamma_{2}} a|\pn u_{2}|^{2}\,\ud x+\int_{\Gamma_{2}} b|u_{2}|^{2}\,\ud x\Big)
\end{equation}
for $h=|\re(\lambda)|^{-1}$ small enough and $|\im(\lambda)|\leq cst$.
\end{lem}
\begin{pr}
In this proof we will keep the same notations as in section~\ref{d15}. We shall extend $w$ in whole $\Omega$ by noting $\tilde{w}=\mathbb{1}_{\Omega_{1}}w_{1}+\mathbb{1}_{\Omega_{2}}w_{2}$ where $\mathbb{1}_{\Omega_{1}}w_{1}$ (resp. $\mathbb{1}_{\Omega_{2}}w_{2}$) is the extension of $w_{1}$ (resp. $w_{2}$) by zero on $\overline{\Omega}_{2}$ (resp. $\overline{\Omega}_{1}$). Note that a such extension is meaningful while $\tilde{w}$ is seen now as a $H^{2}$ function in whole $\Omega$ thanks to the transmission conditions in~\eqref{d7} (See~\cite{D}).

Let $\vfi_{1}$ and $\vfi_{2}$ two weight functions that satisfies the conclusion of the section~\ref{d15}. Let $\chi_{1}=(\chi_{11},\chi_{12})$ and $\chi_{2}=(\chi_{21},\chi_{22})$ two cut-off functions equal to one in $\displaystyle\left(\bigcup_{j=1}^{N_{k}}B(x_{kj},2\epsilon)\right)^{c}$ and supported in $\displaystyle\left(\bigcup_{j=1}^{N_{k}}B(x_{kj},\epsilon)\right)^{c}$ (in order to eliminate the critical points of the phases functions $\vfi_{1}$ and $\vfi_{2}$ (See Figure~\ref{fig2})). Then for $k=1,2$ we obtain from the system~\eqref{d5} the following equations
\begin{equation*}
\left\{
\begin{array}{ll}
\displaystyle\Delta(\chi_{k1}w_{1})=\Psi_{k1}&\text{in }\Omega_{1}
\\
\displaystyle\Delta(\chi_{k2}w_{2})=\Psi_{k2}&\text{in }\Omega_{2}
\\
\chi_{k1}w_{1}=0&\text{on }\Gamma_{1}
\\
\displaystyle \chi_{k2}w_{2}=c_{2}\phi_{a}-ic_{2}\lambda a\,\pn u_{2}-\mathrm{sig}\big(\re(\lambda)\big)\lambda u_{2}&\text{on }\Gamma_{2}
\\
\displaystyle\pn(\chi_{k2}w_{2})=c_{2}\phi_{b}+ic_{2}\lambda b\,u_{2}-\mathrm{sig}\big(\re(\lambda)\big)\lambda\pn u_{2}&\text{on }\Gamma_{2},
\end{array}
\right.
\end{equation*}
where we are noted
\begin{equation}\label{d17}
\Psi_{k}=\left\{\begin{array}{ll}
\displaystyle\Psi_{k1}=[\Delta,\chi_{k1}]w_{1}+\frac{1}{c_{1}}\chi_{k1}\Phi_{1}-\left(\frac{1}{c_{1}h}+\frac{i}{c_{1}}\mathrm{sig}(\re(\lambda))\im(\lambda)\right)(\chi_{k1}w_{1})
\\
\\
\displaystyle\Psi_{k1}=[\Delta,\chi_{k2}]w_{2}+\frac{1}{c_{2}}\chi_{k2}\Phi_{2}-\left(\frac{1}{c_{2}h}+\frac{i}{c_{2}}\mathrm{sig}(\re(\lambda))\im(\lambda)\right)(\chi_{k2}w_{2}),
\end{array}\right.
\end{equation}
and 
\begin{equation*}
\Phi=\left\{
\begin{array}{lcl}
\Phi_{1}&\textrm{in}&\Omega_{1}
\\
\Phi_{2}&\textrm{in}&\Omega_{2}.
\end{array}\right.
\end{equation*}
Applying now Proposition~\ref{d16} to the functions $\chi_{k}\tilde{w}$ and $\Psi_{k}$ with $U=\Omega$ then we obtain for $k=1,2$ that
\begin{equation*}
\begin{split}
h\|\e^{\vfi_{k}/h}\chi_{k}\tilde{w}\|_{L^{2}(U_{k})}^{2}+h^{3}\|\e^{\vfi_{k}/h}\nabla(\chi_{k}\tilde{w})\|_{L^{2}(U_{k})}^{2}\leq
C\big(h^{4}\|\e^{\vfi_{k}/h}\Psi_{k}\|_{L^{2}(U_{k})}^{2}
\\
+h\|\e^{\vfi_{k}/h} w_{2}\|_{L^{2}(\Gamma_{2})}^{2}+h^{3}\|\e^{\vfi_{k}/h}\pn w_{2}\|_{L^{2}(\Gamma_{2})}^{2}\big).
\end{split}
\end{equation*}
Then the expression of $\Psi_{k1}$ and $\Psi_{k2}$ in~\eqref{d17} yields
\begin{equation}\label{d18}
\begin{split}
h\|\e^{\vfi_{k}/h}\chi_{k}w\|_{L^{2}(U_{k})}^{2}+h^{3}\|\e^{\vfi_{k}/h}\nabla(\chi_{k}w)\|_{L^{2}(U_{k})}^{2}\leq
C\big(h^{4}\|\e^{\vfi_{k}/h}\Phi\|_{L^{2}(U_{k})}^{2}
\\
+h^{4}\|\e^{\vfi_{k}/h}[\Delta,\chi_{k}]w\|_{L^{2}(U_{k})}^{2}+h^{4}|\im(\lambda)|^{2}\|\e^{\vfi_{k}/h}\chi_{k}w\|_{L^{2}(U_{k})}^{2}
\\
+h^{3}\|\e^{\vfi_{k}/h}\chi_{k2}w_{2}\|_{L^{2}(U_{k})}^{2}+h\|\e^{\vfi_{k}/h} w_{2}\|_{L^{2}(\Gamma_{2})}^{2}+h^{3}\|\e^{\vfi_{k}/h}\pn w_{2}\|_{L^{2}(\Gamma_{2})}^{2}\big).
\end{split}
\end{equation}
We addition the two last estimates for $k=1,2$ and using the properties of phases $\vfi_{k}<\vfi_{\sigma(k)}$ in $\displaystyle\left(\bigcup_{j=1}^{N_{k}}B(x_{kj},2\epsilon)\right)$ then we can absorb the term $[\Delta,\chi_{k}]w$ at the right hand side of~\eqref{d18} into the left hand side for $h>0$ small. More precisly we obtain
\begin{equation*}
\begin{split}
h\int_{\Omega}\left(\e^{2\vfi_{1}/h}+\e^{2\vfi_{2}/h}\right)|w|^{2}\,\ud x+h^{3}\int_{\Omega}\left(\e^{2\vfi_{1}/h}+\e^{2\vfi_{2}/h}\right)|\nabla w|^{2}\,\ud x
\\
\leq C\Bigg(h^{4}\int_{\Omega}\left(\e^{2\vfi_{1}/h}+\e^{2\vfi_{2}/h}\right)|\Phi|^{2}\,\ud x+h\int_{\Gamma_{2}}\left(\e^{2\vfi_{1}/h}+\e^{2\vfi_{2}/h}\right)|w_{2}|^{2}\,\ud x
\\
+h^{3}\int_{\Gamma_{2}}\left(\e^{2\vfi_{1}/h}+\e^{2\vfi_{2}/h}\right)|\pn w_{2}|^{2}\,\ud x\Bigg).
\end{split}
\end{equation*}
Then by the boundary conditions in~\eqref{d5} we get 
\begin{equation*}
\begin{split}
\int_{\Omega}|w|^{2}\,\ud x+\int_{\Omega}|\nabla w|^{2}\,\ud x\leq C\e^{C/h}\Bigg(\int_{\Omega}|\Phi|^{2}\,\ud x+\int_{\Gamma_{2}}|u_{2}|^{2}\,\ud x+\int_{\Gamma_{2}}|b\,u_{2}|^{2}\,\ud x
\\
+\int_{\Gamma_{2}}|\pn u_{2}|^{2}+\int_{\Gamma_{2}}|a\,\pn u_{2}|^{2}\,\ud x\Bigg).
\end{split}
\end{equation*}
And this yields from the assumption~\eqref{2} that
\begin{equation}\label{d12}
\|w\|_{L^{2}(\Omega)}^{2}\leq C\e^{C/h}\Big(\|\Phi\|_{L^{2}(\Omega)}^{2}+\|\phi_{a}\|_{L^{2}(\Gamma_{2})}^{2}+\|\phi_{b}\|_{L^{2}(\Gamma_{2})}^{2}+\int_{\Gamma_{2}} a|\pn u_{2}|^{2}\,\ud x+\int_{\Gamma_{2}} b|u_{2}|^{2}\,\ud x\Big).
\end{equation}
Observing by Green's formula and the expression of $w$ in~\eqref{d7} that
\begin{equation}\label{d19}
\|w\|_{L^{2}(\Omega)}^{2}=\alpha\|\Delta u\|_{L^{2}(\Omega)}^{2}+|\lambda|^{2}\|u\|_{L^{2}(\Omega)}^{2}+2\alpha|\re(\lambda)|.\|\nabla u\|_{L^{2}(\Omega)}^{2}\geq C\|u\|_{X}^{2}.
\end{equation}
And this completes the proof by the combination of~\eqref{d12} and~\eqref{d19}.
\end{pr}

From~\eqref{d9} and~\eqref{d10} we obtain
\begin{equation*}
\begin{split}
\|u\|_{X}^{2}\leq C\e^{C/h}\Big(\|\Phi\|_{H}^{2}+\|\phi_{a}\|_{L^{2}(\Gamma_{2})}^{2}+\|\phi_{b}\|_{L^{2}(\Gamma_{2})}^{2}+\|\Phi\|_{H}\|u\|_{H}
\\
+|\re(\lambda)\im(\lambda)|^{2}\|u\|_{H}^{2}+\|f\|_{2}\|u\|_{2}\Big),
\end{split}
\end{equation*}
then by the expression of $\phi_{a}$ and $\phi_{b}$ in~\eqref{d4} we have
\begin{equation}\label{d11}
\|u\|_{X}^{2}\leq C\e^{C/h}\left(\|\Phi\|_{H}^{2}+\|f_{2}\|_{2}^{2}+\|\Phi\|_{H}\|u\|_{H}+|\im(\lambda)|^{2}\|u\|_{H}^{2}+\|f\|_{2}\|u\|_{2}\right).
\end{equation}
Then estimate~\eqref{d11} and Proposition~\ref{7} for $|\im(\lambda)|\leq \frac{1}{\sqrt{2C}}\e^{-C/h}$ give us
$$\|u\|_{X}^{2}\leq C\e^{C/h}\left(\|\Phi\|_{H}^{2}+\|f_{2}\|_{X}^{2}\right).$$
Using the expression of $\Phi$ in~\eqref{d4} we obtain
\begin{equation}\label{d13}
\|u\|_{X}\leq C\e^{C/h}\left(\|F\|_{X}+\|G\|_{H}\right).
\end{equation}
We thus obtain form the first equation of~\eqref{d7} and~\eqref{d13} that
\begin{equation}\label{d14}
\|v\|_{H}\leq |\lambda|\,\|u\|_{H}+\|F\|_{H}\leq C\e^{C/h}\left(\|F\|_{X}+\|G\|_{H}\right),
\end{equation}
and hence~\eqref{d13} and~\eqref{d14} give
$$\|(u,v)\|_{\mathcal{H}}\leq C\e^{C|\re(\lambda)|}\|(i\mathcal{A}+\lambda\id)(u,v)\|_{\mathcal{H}}^{2}.$$
Then $(i\mathcal{A}+\lambda)$ is injective then bijective in $\mathcal{D}(\mathcal{A})$ and we have
$$\|(i\mathcal{A}+\lambda\id)^{-1}\|_{\mathscr{L}(\mathcal{H},\mathcal{H})}\leq C\e^{C|\re(\lambda)|}$$
for $\lambda\in\{z\in\mathbb{C};\;|\im(z)|<C_{1}\e^{-C_{2}|\re(z)|},\,|z|>C_{3}\}$ and this complete the proof of Theorem~\ref{9}.
\nocite{*}
\bibliographystyle{alpha}
\bibliography{biblio1}
\addcontentsline{toc}{section}{References}
\end{document}